\documentclass[12pt]{elsarticle}
\usepackage{hyperref}
\usepackage[textwidth=16cm, textheight=24cm, left=2.5cm, right=2.5cm, top=3cm, bottom=3cm]{geometry}
\usepackage{latexsym}
\usepackage{enumerate}
\usepackage{amsmath, amsthm}
\usepackage{amssymb}
\usepackage{graphics}
\usepackage{graphicx}
\usepackage{subfigure}
\usepackage{indentfirst}
\usepackage{mathrsfs}
\usepackage{tikz}
\usepackage{longtable}
\usepackage{multirow}

\usetikzlibrary{decorations.pathreplacing}

\theoremstyle{plain}

\theoremstyle{definition}

\numberwithin{equation}{section}

\makeatletter
\def\ps@pprintTitle{%
  \let\@oddhead\@empty
  \let\@evenhead\@empty
  \def\@oddfoot{\reset@font\hfil\thepage\hfil}
  \let\@evenfoot\@oddfoot
}
\makeatother

\begin{document}
\begin{frontmatter}
\title{Some Spectral Properties and Characterizations of Connected Odd-bipartite Uniform Hypergraphs \tnoteref{title1}}
\tnotetext[title1]{Research supported by the National Science Foundation of China No.11231004,11271288 and 11301340}
\author[a]{Jia-Yu Shao}
\ead{jyshao@tongji.edu.cn}
\author[a]{Hai-Ying Shan}
\ead{shan\_haiying@tongji.edu.cn}
\author[b]{Bao-feng Wu}
\ead{baufern@aliyun.com}

\address[a]{Department of Mathematics,  Tongji University,  Shanghai  200092,  China}
\address[b]{College of Science, University of Shanghai for Science and Technology, Shanghai, 200093, China}

\date{}
\begin{abstract}
A $k$-uniform hypergraph $G=(V,E)$ is called odd-bipartite ([5]), if $k$ is even and there exists some proper subset $V_1$ of $V$ such that each edge of $G$ contains odd number of vertices in $V_1$. Odd-bipartite hypergraphs are generalizations of the ordinary bipartite graphs. We study the spectral properties of the connected odd-bipartite hypergraphs. We prove that the Laplacian H-spectrum and signless Laplacian H-spectrum of a connected $k$-uniform hypergraph $G$ are equal if and only if $k$ is even and $G$ is odd-bipartite. We further give several spectral characterizations of the connected odd-bipartite hypergraphs. We also give a characterization for a connected $k$-uniform hypergraph whose Laplacian spectral radius and signless Laplacian spectral radius are equal, thus provide an answer to a question raised in [9]. By showing that the Cartesian product $G\Box H$ of two odd-bipartite $k$-uniform hypergraphs is still odd-bipartite, we determine that the Laplacian spectral radius of $G\Box H$ is the sum of the Laplacian spectral radii of $G$ and $H$, when $G$ and $H$ are both connected odd-bipartite.

\end{abstract}

\end{frontmatter}

\vskip 0.68cm

\section{Introduction}

In recent years, the study of tensors and the spectra of tensors (and hypergraphs) with their various applications has attracted extensive attention and interest, since the work of L.Qi ([8]) and L.H.Lim ([7]) in 2005.

\vskip 0.1cm

As was in [8], an order $m$ dimension $n$ tensor $\mathbb {A}=(a_{i_1i_2\cdots i_m})_{1\le i_j\le n \ (j=1,\cdots ,m)}$ over the complex field $\mathbb {C}$ is a multidimensional array with all entries $a_{i_1i_2\cdots i_m}\in \mathbb {C} \  \ (i_1,\cdots ,i_m\in [n]=\{1,\cdots ,n\})$.

In this paper, we only consider real tensors. We first need the following definition of a general product of two tensors.

\vskip 0.18cm

\noindent {\bf Definition 1.1}([10]): Let $\mathbb {A}$ (and $\mathbb {B}$) be an order $m\ge 2$ (and order $k\ge 1$), dimension $n$ tensor, respectively. Define the product $\mathbb {A}\cdot \mathbb {B}$ (or simply $\mathbb {A} \mathbb {B}$) to be the following tensor $\mathbb {C}$ of order $(m-1)(k-1)+1$ and dimension $n$:
$$c_{i\alpha_1\cdots \alpha_{m-1}}=\sum_{i_2,\cdots ,i_m=1}^n a_{ii_2\cdots i_m}b_{i_2\alpha_1}\cdots b_{i_m\alpha_{m-1}} \quad  (i\in [n], \  \alpha_1, \cdots ,\alpha_{m-1}\in [n]^{k-1}) $$
It is proved in [10] that this product of tensors satisfies the associative law.

\vskip 0.18cm

In particular, when $P$ and $Q$ are both matrices and $\mathbb {B}=P\mathbb {A}Q$, we have:

$$b_{i_1\cdots i_m}=\sum_{j_1,\cdots ,j_m=1}^n a_{j_1\cdots j_m}p_{i_1j_1}q_{j_2i_2}\cdots q_{j_mi_m}  \eqno {(1.1)} $$

\vskip 0.18cm

Let $\mathbb {A}$  be an order $m$  dimension $n$ tensor, let $x=(x_1,\cdots ,x_n)^T\in \mathbb {C}^n$ be a column vector of dimension $n$. Then their product $\mathbb {A}x$ is a vector in $\mathbb {C}^n$ whose $i$th component is as the following:
$$(\mathbb {A}x)_i=\sum_{i_2,\cdots ,i_m=1}^n a_{ii_2\cdots i_m}x_{i_2}\cdots x_{i_m} \eqno {(1.2)}$$
Let $x^{[r]}=(x_1^r,\cdots ,x_n^r)^T$. Then ([1,8]) a number $\lambda \in \mathbb {C}$ is called an eigenvalue of the tensor $\mathbb {A}$ if there exists a nonzero vector $x\in \mathbb {C}^n$ such that
$$\mathbb {A}x=\lambda x^{[m-1]} \eqno {(1.3)}$$
and in this case, $x$ is called an eigenvector of $\mathbb {A}$ corresponding to the eigenvalue $\lambda$.

\vskip 0.18cm

An eigenvalue of $\mathbb {A}$ is called an H-eigenvalue ([8]), if there exists a real eigenvector corresponding to it.

\vskip 0.18cm

The maximal absolute value of the eigenvalues of $\mathbb {A}$ is called the spectral radius of $\mathbb {A}$, denoted by $\rho (\mathbb {A})$. The largest H-eigenvalue of a real symmetric tensor $\mathbb {A}$ is denoted by $\lambda (\mathbb {A})$.

\vskip 0.18cm

The equation (1.3) can also be written in the following equivalent form:
$$(\lambda \mathbb {I}- \mathbb {A})\cdot x=0  \eqno {(1.4)}$$
where $\mathbb {I}$ is the unit tensor of order $m$ and dimension $n$.

\vskip 0.18cm

In order to define the characteristic polynomials and spectra of tensors, we first need to define the determinants of tensors. Originally the determinants of tensors were defined as the resultants of some corresponding system of homogeneous equations on $n$ variables. Here we give the following equivalent definition of the determinants of tensors by using some properties of the determinants proved in [4].

\vskip 0.38cm

\noindent {\bf Definition 1.2} ([4]): Let $\mathbb {A}$ be an order $m$ dimension $n$ tensor with $m\ge 2$. Then its determinant $det(\mathbb {A})$ is defined to be the unique polynomial on the entries of $\mathbb {A}$ satisfying the following three conditions:
\vskip 0.1cm

\noindent (1) $det(\mathbb {A})=0$ if and only if the system of homogeneous equations $\mathbb {A}x=0$ has a nonzero solution.

\vskip 0.1cm

\noindent (2) $det(\mathbb {A})=1$, when $\mathbb {A}=\mathbb {I}$ is the unit tensor.

\vskip 0.1cm

\noindent (3) $det(\mathbb {A})$ is an irreducible polynomial on the entries of $\mathbb {A}$, when the entries of $\mathbb {A}$ are viewed as dsitinct independent variables.

\vskip 0.38cm

\noindent {\bf Definition 1.3}: Let $\mathbb {A}$ be an order $m\ge 2$ dimension $n$ tensor. Then the characteristic polynomial of $\mathbb {A}$, denoted by $\phi_{\mathbb {A}}(\lambda )$, is the determinant $det(\lambda \mathbb {I}- \mathbb {A})$.

\vskip 0.18cm

From the above definitions, it is easy to see that $\lambda$ is an eigenvalue of $\mathbb {A}$ if and only if it is a root of the characteristic polynomial of $\mathbb {A}$.

\vskip 0.28cm

\noindent {\bf Definition 1.4}: Let $\mathbb {A}$ be an order $m\ge 2$ dimension $n$ tensor. Then the (multi)-set of roots of the characteristic polynomial of $\mathbb {A}$ (counting multiplicities) is called the spectrum of $\mathbb {A}$, denoted by $Spec(\mathbb {A})$.

The H-spectrum of a real tensor $\mathbb {A}$, denoted by $Hspec (\mathbb {A})$, is defined to be the set of distinct H-eigenvalues of $\mathbb {A}$. Namely,
$$Hspec (\mathbb {A})=\{\lambda \in \mathbb {R} \ | \ \lambda  \ \mbox {is an H-eigenvalue of $\mathbb {A}$} \ \}$$

\vskip 0.18cm

\noindent {\bf Definition 1.5} ([10, 11]): Let $\mathbb {A}$ and $\mathbb {B}$ be two order $m\ge 2$  dimension $n$ tensors. Suppose that there exists a nonsingular diagonal matrix $D$ of order $n$ such that $\mathbb {B} = D^{-(m-1)}\mathbb {A}D  $,
then $\mathbb {A}$ and $\mathbb {B}$ are called diagonal similar.

\vskip 0.28cm

It is proved in [10, Theorem 2.1] that similar tensors (thus diagonal similar tensors) have the same characteristic polynomials, and thus have the same spectra (just as in the case of matrices).

\vskip 0.28cm

In [3], Friedland et al. defined the weak irreducibility of nonnegative tensors. It was proved in [3] and [11] that a $k$-hypergraph $G$ is connected if and only if its adjacency tensor $\mathbb {A}(G)$ is weakly irreducible. They further generalized the results of the well known Perron-Frobenius Theorem from nonnegative irreducible matrices to nonnegative weakly irreducible tensors.

\vskip 0.28cm

\noindent {\bf Lemma 1.1}([3,11]): Let $\mathbb {A}$ be a nonnegative tensor. Then

\vskip 0.1cm

\noindent {\bf (1)}. $\rho (\mathbb {A})$ is an H-eigenvalue of $\mathbb {A}$ with a nonnegative eigenvector. Furthermore, if $\mathbb {A}$ is weakly irreducible, then $\rho (\mathbb {A})$ has a positive eigenvector.

\vskip 0.1cm

\noindent {\bf (2)}. If $\lambda$ is an eigenvalue of $\mathbb {A}$ with a positive eigenvector, then $\lambda=\rho (\mathbb {A})$.

\vskip 0.28cm

\noindent {\bf Lemma 1.2}([11]): Let $\mathbb {A}$ and $\mathbb {B}$ be two order  $k$ dimension $n$ tensors with $|\mathbb {B}|\le \mathbb {A}$. Then

\vskip 0.1cm

\noindent {\bf (1)}. $\rho (\mathbb {B}) \le \rho (\mathbb {A})$.

\vskip 0.1cm

\noindent {\bf (2)}. Furthermore, if $\mathbb {A}$ is weakly irreducible and $\rho (\mathbb {B}) = \rho (\mathbb {A})$, where $\lambda =\rho (\mathbb {A})e^{i\varphi}$ is an eigenvalue of $\mathbb {B}$ with an eigenvector $y$, then

\vskip 0.1cm

\noindent {\bf (i)}. All the components of $y$ are nonzero.

\vskip 0.1cm

\noindent {\bf (ii)}. Let $U=diag (y_1/|y_1|,\cdots, y_n/|y_n|)$ be a nonsingular diagonal matrix, then we have
$$\mathbb {B} = e^{i\varphi}U^{-(k-1)}\mathbb {A}U  $$

\vskip 0.38cm

A $k$-uniform hypergraph (or simply a $k$-hypergraph) $G$ is a hypergraph each of whose edges contains exactly $k$ vertices. The adjacency tensor of $G$ (under certain ordering of vertices) is the order $k$ dimension $n$ tensor $\mathbb {A}=\mathbb {A}(G)$ with the following entries ([2]):
$$a_{i_1i_2\cdots i_k}=\left \{\begin{array}{cc}
                               \frac {1} {(k-1)!} &   \mbox{if }  \ \{i_1,i_2,\cdots ,i_k\} \in E(G) \\
                               0 &   \mbox{otherwise}\\
                            \end{array}
                          \right.$$
The characteristic polynomial and spectrum of a uniform hypergraph $G$ are those of its adjacency tensor $\mathbb {A}$.

\vskip 0.18cm

Let $\mathbb {D}=\mathbb {D}(G)$ be the degree diagonal tensor of $G$ (its $i$-th diagonal element is the degree of the vertex $i$), then the tensor $\mathbb {L}=\mathbb {D}-\mathbb {A}$ is called the Laplacian tensor of $G$, and $\mathbb {Q}=\mathbb {D}+\mathbb {A}$ is called the signless Laplacian tensor of $G$. The Laplacian spectrum and signless Laplacian spectrum of $G$ are defined to be the spectrum of $\mathbb {L}$ and $\mathbb {Q}$, respectively.

\vskip 0.28cm

\noindent {\bf Definition 1.6}([5]): A $k$-hypergraph $G=(V,E)$ is called odd-bipartite, if $k$ is even and there exists some proper subset $V_1$ of $V$ such that each edge of $G$ contains exactly odd number of vertices in $V_1$.

\vskip 0.18cm

It is easy to see that in the case $k=2$, a $k$-hypergraph $G$ is odd-bipartite if and only if $G$ is an ordinary bipartite graph. Thus the concept of odd-bipartite hypergraphs is a generalization of that of (ordinary) bipartite graphs.

\vskip 0.18cm

In [5, Proposition 3.2], Hu and Qi proved that if $G$ is a $k$-uniform hm-bipartite hypergraph with $k$ even (which is a special subclass of the class of odd-bipartite $k$-hypergraphs), then its Laplacian spectrum and signless Laplacian spectrum are equal. They also ask the question about whether the converse of this result is true or not. In Section 2, we study the following analogy question for Laplacian H-spectrum and signless Laplacian H-spectrum:

\vskip 0.28cm

\noindent {\bf Question 1}: Characterize those connected $k$-hypergraphs whose Laplacian H-spectra and signless Laplacian H-spectra are equal.

\vskip 0.28cm

We will show in Section 2 that in the case when $G$ is connected, the Laplacian H-spectrum and signless Laplacian H-spectrum of $G$ are equal if and only if $k$ is even and $G$ is odd-bipartite. This result is a generalization of the corresponding spectral characterization of the ordinary bipartite graphs. In Section 2, we also show that the equality of the Laplacian H-spectrum and signless Laplacian H-spectrum implies the equality of the Laplacian spectrum and signless Laplacian spectrum, and give several further spectral characterizations of the connected odd-bipartite hypergraphs in Theorems 2.1-2.5. These characterizations show that some of the hypergraph structures can be well described by the spectra of the hypergraphs, including the adjacency spectra, the Laplacian spectra and signless Laplacian spectra of the hypergraphs, just as in the case of ordinary graphs ($k=2$).

\vskip 0.28cm

In [9], it was mentioned that ``It is a research topic to identify the conditions under which $\rho (\mathbb {L}) = \rho (\mathbb {Q})$". In Theorem 2.4, we use the properties of the weakly irreducible tensors and the diagonal similar tensors to show that when $G$ is a connected $k$-hypergraph, then $\rho (\mathbb {L}) = \rho (\mathbb {Q})$ if and only if $Spec(\mathbb {L})= Spec (\mathbb {Q})$.

\vskip 0.18cm

In Section 3, we will give the expressions of the Laplacian tensor and signless Laplacian tensor of the Cartesian product $G\Box H$ of $k$-hypergraphs by using the direct product of tensors defined in [10]. By showing that the Cartesian product of two odd-bipartite $k$-hypergraphs is still odd-bipartite together with the formula $(\mathbb {A}\otimes \mathbb {B})(\mathbb {C}\otimes \mathbb {D})=(\mathbb {AC})\otimes (\mathbb {BD})$  for the direct products and general products of tensors, we determine that the Laplacian spectral radius (and the largest Laplacian H-eigenvalue) of $G\Box H$ is the sum of the Laplacian spectral radii (and the largest Laplacian H-eigenvalues) of $G$ and $H$, when $G$ and $H$ are both connected odd-bipartite.

\vskip 1.88cm

\section{The relation between the Laplacian spectra and signless Laplacian spectra of odd-bipartite hypergraphs}

\vskip 0.18cm

The following theorem gives an algebraic feature of the odd-bipartite hypergraphs by using the tensor product defined in Section 1.

\vskip 0.18cm

\noindent {\bf Theorem 2.1}: Let $G$ be a connected $k$-uniform hypergraph with $n$ vertices, $\mathbb {A}$ and $\mathbb {D}$ be the adjacency tensor and degree diagonal tensor of $G$, $\mathbb {L}=\mathbb {D}-\mathbb {A}$ and $\mathbb {Q}=\mathbb {D}+\mathbb {A}$ be the Laplacian and signless Laplacian tensor of $G$, respectively. Then the following three conditions are equivalent:

\vskip 0.18cm

\noindent (1). There exists some diagonal matrix $P$ of order $n$ with all the diagonal entries $\pm 1$ and $P\ne -I_n$ such that $\mathbb {L}=P^{-(k-1)}\mathbb {Q}P$.

\vskip 0.18cm

\noindent (2). There exists some diagonal matrix $P$ of order $n$ with all the diagonal entries $\pm 1$ and $P\ne -I_n$ such that $\mathbb {A}= -P^{-(k-1)}\mathbb {A}P$.

\vskip 0.18cm

\noindent (3). $k$ is even and $G$ is odd-bipartite.

\vskip 0.18cm

\noindent {\bf Proof}. (1) $\Longleftrightarrow$ (2): For the diagonal tensor $\mathbb {D}$ of order $k$ and diagonal matrix $P$ of order $n$ with all the diagonal entries $\pm 1$, it is easy to see that we always have $\mathbb {D}=P^{-(k-1)}\mathbb {D}P$. Thus we have:

$$\mathbb {L}=P^{-(k-1)}\mathbb {Q}P \Longleftrightarrow \mathbb {D}-\mathbb {A}=P^{-(k-1)}(\mathbb {D}+\mathbb {A})P\Longleftrightarrow \mathbb {A}=-P^{-(k-1)}\mathbb {A}P$$

\vskip 0.1cm

(2) $\Longrightarrow$ (3):  Since $G$ is nontrivial, we have $\mathbb {A}\ne O$. So $P$ is also not the identity matrix, for otherwise we would have $\mathbb {A}=-\mathbb {A}$, thus $\mathbb {A}=O$, a contradiction. Let
$$V_1=\{i\in [n] \ | \ \mbox {the $i$-th diagonal entry of $P$ is $-1$} \}$$
Then $V_1$ is a proper subset of [n].

\vskip 0.1cm

Now by $\mathbb {A}= -P^{-(k-1)}\mathbb {A}P$ and using (1.1) we have
$$-a_{i_1i_2\cdots i_k}=p_{i_1i_1}^{-(k-1)}p_{i_2i_2}\cdots p_{i_ki_k}a_{i_1i_2\cdots i_k} \qquad (\forall \  1\le i_1,i_2,\cdots, i_k \le n) \eqno {(2.1)}$$

Let $e=\{i_1,i_2,\cdots, i_k\}$ be an edge of $G$. Then $a_{i_1i_2\cdots i_k}\ne 0$ and so from (2.1) we have

$$- p_{i_1i_1}^k=p_{i_1i_1}p_{i_2i_2}\cdots p_{i_ki_k} \qquad ( \ \forall \  \{i_1,i_2,\cdots, i_k\}\in E(G) \ ) \eqno {(2.2)}$$
which implies that
$$p_{ii}^k=p_{jj}^k \qquad \mbox {(if $i$ and $j$ are adjacent in $G$)} \eqno {(2.3)}$$

Now if $k$ is odd, then from (2.3) we would have
$$p_{ii}=p_{jj} \qquad \mbox {(if $i$ and $j$ are adjacent in $G$)} $$
From this we deduce that all the diagonal entries of $P$ are equal since $G$ is connected, which would imply that either $P$ or $-P$ is the identity matrix, a contradiction. From this we conclude that $k$ is even.

Now suppose that $e=\{i_1,i_2,\cdots, i_k\}$ is an edge of $G$ with $|e\cap V_1|=r$. Then $a_{i_1i_2\cdots i_k}\ne 0$ and so from (2.2) we further have $-1=(-1)^r$ since $k$ is even, thus $r$ is odd. From this we conclude that $G$ is odd-bipartite.

\vskip 0.1cm

(3) $\Longrightarrow$ (2): Suppose that $k$ is even and $G$ is odd-bipartite, then there exists some proper subset $V_1$ of [n] such that every edge of $G$ intersects $V_1$ with exactly an odd number of vertices.

Now take $P$ to be the diagonal matrix of order $n$ with all the diagonal entries $\pm 1$ such that $p_{ii}=-1$ if and only if $i\in V_1$. Then $P\ne -I_n$, and we can check that (2.1) holds since $k$ is even, which means that $-\mathbb {A}=P^{-(k-1)}\mathbb {A}P$.
\qed

\vskip 0.18cm

Using this algebraic characterization of odd-bipartite hypergraphs, we are now able to obtain the following characterization for a connected $k$-uniform hypergraph whose Laplacian H-spectrum and signless Laplacian H-spectrum are equal. Thus obtain an answer to the Question 1 in Section 1.

\vskip 0.38cm

\noindent {\bf Theorem 2.2}: Let $G$ be a connected $k$-uniform hypergraph with $n$ vertices. Let $\mathbb {A}$ and $\mathbb {D}$ be the adjacency tensor and degree diagonal tensor of $G$, and $\mathbb {L}= \mathbb {D}-\mathbb {A}$ and $\mathbb {Q}= \mathbb {D}+\mathbb {A}$ be the Laplacian and signless Laplacian tensor of $G$, respectively. Then the following three conditions are equivalent:

\vskip 0.1cm

\noindent (1). $k$ is even and $G$ is odd-bipartite.

\vskip 0.1cm

\noindent (2).  $Spec(\mathbb {L})= Spec (\mathbb {Q})$ and $Hspec(\mathbb {L})= Hspec (\mathbb {Q})$.

\vskip 0.1cm

\noindent (3).  $Hspec(\mathbb {L})= Hspec (\mathbb {Q})$.

\vskip 0.18cm

\noindent {\bf Proof.} (1)$\Longrightarrow $(2): By Theorem 2.1, we see that (1) implies that there exists  some diagonal matrix $P$ of order $n$ with all the diagonal entries $\pm 1$ and $P\ne -I_n$, such that
$$\mathbb {L}=P^{-(k-1)}\mathbb {Q}P $$
Thus $\mathbb {L}$ and $\mathbb {Q}$ are diagonal similar (in the sense of Definition 1.5), so by [10] we have $Spec(\mathbb {L})= Spec (\mathbb {Q})$.

\vskip 0.1cm

On the other hand, Let $y=Px$. Then by the relation $\mathbb {L}=P^{-(k-1)}\mathbb {Q}P $ and the associativity of the product of tensors we can check that
$$\mathbb {L}x=\lambda x^{[k-1]}\Longleftrightarrow P^{-(k-1)}\mathbb {Q}Px=\lambda x^{[k-1]}\Longleftrightarrow \mathbb {Q}y=P^{(k-1)}\lambda x^{[k-1]}=\lambda (Px)^{[k-1]}$$

\quad $\Longleftrightarrow \mathbb {Q}y=\lambda y^{[k-1]}$

\noindent Since $P$ is a real nonsingular matrix, the above relation implies that $\lambda$ is an H-eigenvalue of $\mathbb {L}$ if and only if it is an H-eigenvalue of $\mathbb {Q}$. Thus we also have  $Hspec(\mathbb {L})= Hspec (\mathbb {Q})$.

\vskip 0.18cm

(2)$\Longrightarrow $(3): This is obvious.

\vskip 0.18cm

(3)$\Longrightarrow $(1): Since $G$ is connected,  $\mathbb {Q}$ is nonnegative weakly irreducible. So by Lemma 1.1 we know that $\rho (\mathbb {Q})$ is an H-eigenvalue of $\mathbb {Q}$, and thus an H-eigenvalue of $\mathbb {L}$ by condition (3).

From this we also have $\rho (\mathbb {L})=\rho (\mathbb {Q})$, since $|\mathbb {L}|= \mathbb {Q}$ implying $\rho (\mathbb {L})\le \rho (\mathbb {Q})$ by (1) of Lemma 1.2.

Let $y=(y_1,y_2,\cdots,y_n)^T$ be a real eigenvector of $\mathbb {L}$ corresponding to the H-eigenvalue $\rho (\mathbb {L})= \rho (\mathbb {Q})$. Then all the components of $y$ are nonzero by Lemma 1.2. Without loss of generality we may assume that $y_1> 0$. Let $P=diag (y_1/|y_1|,\cdots, y_n/|y_n|)$, then $P$ is a diagonal matrix of order $n$ with all the diagonal entries $\pm 1$ and $P\ne -I_n$, since $y$ is real and $y_1> 0$.

Note that $\mathbb {Q}$ is nonnegative weakly irreducible and $|\mathbb {L}|= \mathbb {Q}$, so by Lemma 1.2 we have $\mathbb {L}=P^{-(k-1)}\mathbb {Q}P $. Thus by Theorem 2.1 we conclude that $k$ is even and $G$ is odd bipartite.
\qed

\vskip 0.38cm

Similarly, we can obtain the following characterizations of the (connected) odd-bipartite $k$-hypergraphs in terms of the symmetry of their adjacency spectra and adjacency H-spectra.

\vskip 0.18cm

\noindent {\bf Theorem 2.3}: Let $G$ be a connected $k$-uniform hypergraph with $n$ vertices, and $\mathbb {A}$ be the adjacency tensor of $G$. Then the following three conditions are equivalent:

\vskip 0.1cm

\noindent (1). $k$ is even and $G$ is odd-bipartite.

\vskip 0.1cm

\noindent (2).  $Spec(\mathbb {A})= -Spec (\mathbb {A})$ and $Hspec(\mathbb {A})= -Hspec (\mathbb {A})$. Namely, both $Spec(\mathbb {A})$ and $Hspec(\mathbb {A})$ are symmetric with respect to the origin.

\vskip 0.1cm

\noindent (3).  $Hspec(\mathbb {A})= -Hspec (\mathbb {A})$.

\vskip 0.18cm

\noindent {\bf Proof.} The proof is similar to that of Theorem 2.2. The only difference is that now we need to use the part (2)$\Longleftrightarrow $(3) of Theorem 2.1 instead of the part (1)$\Longleftrightarrow $(3).
\qed

\vskip 0.38cm

Theorems 2.1-2.3 can serve as the examples to show how the structures of hypergraphs (e.g., the odd-bipartite property of hypergraphs) can be described and determined by the spectral properties (in particular, the Laplacian and signless Laplacian spectral properties) of the hypergraphs. This may also be viewed as one of the advantages for the study of the Laplacian and signless Laplacian spectra of hypergraphs.

\vskip 0.38cm

In [9], it was mentioned that ``It is a research topic to identify the conditions under which $\rho (\mathbb {L}) = \rho (\mathbb {Q})$". In the following Theorem 2.4, we will show that when $G$ is a connected $k$-hypergraph, then $\rho (\mathbb {L}) = \rho (\mathbb {Q})$ if and only if $Spec(\mathbb {L})= Spec (\mathbb {Q})$.

\vskip 0.28cm

\noindent {\bf Theorem 2.4}: Let $G$ be a connected $k$-hypergraph, $\mathbb {L}= \mathbb {D}-\mathbb {A}$ and $\mathbb {Q}= \mathbb {D}+\mathbb {A}$ be the Laplacian and signless Laplacian tensors of $G$, respectively. Then $\rho (\mathbb {L}) = \rho (\mathbb {Q})$ if and only if $Spec(\mathbb {L})= Spec (\mathbb {Q})$.

\vskip 0.18cm

\noindent {\bf Proof.} The sufficiency part is obvious. We now prove the necessity.

Suppose $\rho (\mathbb {L}) = \rho (\mathbb {Q})$, where $\lambda =\rho (\mathbb {Q})e^{i\varphi}$ is an eigenvalue of $\mathbb {L}$. Since $|\mathbb {L}|= \mathbb {Q}$ and $\mathbb {Q}$ is a nonnegative weakly irreducible tensor, by Lemma 1.2 we know that there exists a nonsingular diagonal matrix $U$ such that $\mathbb {L}= e^{i\varphi}U^{-(k-1)}\mathbb {Q}U$, namely
$$\mathbb {D}-\mathbb {A} = e^{i\varphi}U^{-(k-1)}(\mathbb {D}+\mathbb {A})U  \eqno {(2.4)} $$
Also, it can be verified that $U^{-(k-1)}\mathbb {D}U= \mathbb {D}$, since the matrix $U$ and the tensor $\mathbb {D}$ are both diagonal. So by comparing the diagonal entries of the both sides of (2.4), we have
$$\mathbb {D} = e^{i\varphi}U^{-(k-1)}\mathbb {D}U = e^{i\varphi}\mathbb {D}  $$
Thus we have $e^{i\varphi}=1$, and so $\mathbb {L}= U^{-(k-1)}\mathbb {Q}U$, which means that $\mathbb {L}$ and $\mathbb {Q}$ are diagonal similar. So by [10, Theorem 2.3] we have $Spec(\mathbb {L})= Spec (\mathbb {Q})$.
\qed

\vskip 0.38cm

As an application of Theorem 2.2, we will give another spectral characterization of connected odd-bipartite $k$-hypergraphs in the following Theorem 2.5. This result can be proved by using [5]. Here we will give a different proof which uses Theorem 2.2. We first have the following elementary lemma.

\vskip 0.18cm

\noindent {\bf Lemma 2.1}: Let $k\ge 2$ be an even integer, $x_1,\cdots, x_k$ be real numbers. Then we have
$$x_1^k+\cdots + x_k^k\pm kx_1\cdots x_k\ge 0$$

\noindent {\bf Proof.} By the fundamental inequality $\frac {1}{k}(a_1+\cdots+a_k)\ge (a_1\cdots a_k)^{\frac {1}{k}} \ $ (for $a_1,\cdots, a_k\ge 0$), we have
$$\frac {1}{k}(x_1^k+\cdots + x_k^k)=\frac {1}{k}(|x_1|^k+\cdots + |x_k|^k)\ge |x_1\cdots x_k| \ge  \pm x_1\cdots x_k$$
\qed

\vskip 0.18cm

\noindent {\bf Theorem 2.5}: Let $G$ be a connected $k$-hypergraph, $\mathbb {Q}= \mathbb {D}+\mathbb {A}$ be the signless Laplacian tensor of $G$. Then the following two conditions are equivalent:

\vskip 0.1cm

\noindent (1). $k$ is even and $G$ is odd-bipartite.

\vskip 0.1cm

\noindent (2).  0 is an H-eigenvalue of $\mathbb {Q}$.

\vskip 0.18cm

\noindent {\bf Proof.}  (1)$\Longrightarrow $(2): First, it is easy to see that 0 is an H-eigenvalue of the Laplacian tensor $\mathbb {L}= \mathbb {D}-\mathbb {A}$, with the all 1 vector as its (real) eigenvector.

Second, from Theorem 2.2 we know that (1)$\Longrightarrow Hspec(\mathbb {L})= Hspec (\mathbb {Q})$. Thus (1) also implies that  0 is an H-eigenvalue of $\mathbb {Q}$.

\vskip 0.1cm

(2)$\Longrightarrow $(1): By [5, Proposition 4.1] we know that (2)$\Longrightarrow k$ is even.

Now let $x\in R^n$ be a real eigenvector of $\mathbb {Q}$ corresponding to the H-eigenvalue 0. Then we have $\mathbb {Q}x=0$. Thus we have
$$0=x^T(\mathbb {Q}x)=\sum_{i_1,\cdots,i_k=1}^n q_{i_1\cdots i_k}x_{i_1}\cdots x_{i_k}$$
On the other hand, it is not difficult to calculate that

$$\begin{aligned}
\displaystyle
0&=\sum_{i_1,\cdots,i_k=1}^n q_{i_1\cdots i_k}x_{i_1}\cdots x_{i_k} =\sum_{i_1,\cdots,i_k=1}^n d_{i_1\cdots i_k}x_{i_1}\cdots x_{i_k}+ \sum_{i_1,\cdots,i_k=1}^n a_{i_1\cdots i_k}x_{i_1}\cdots x_{i_k} \\
&=\sum_{\{i_1,\cdots,i_k\} \in E(G)}(x_{i_1}^k+\cdots + x_{i_k}^k + kx_{i_1}\cdots x_{i_k})\\
\end{aligned} $$
So by Lemma 2.1 we have
$$x_{i_1}^k+\cdots + x_{i_k}^k + kx_{i_1}\cdots x_{i_k}=0 \qquad (\forall \ \{i_1,\cdots,i_k\} \in E(G))\eqno {(2.5)}$$

From (2.5) we see that, if $\{i_1,\cdots,i_k\} \in E(G)$ with $x_{i_1}=0$, then we would have
$x_{i_2}=\cdots =x_{i_k}=0$. This means that if some component of $x$ is zero, then by the connectivity of $G$ all the components of $x$ should be zero. This contradicts the fact that $x$ is an eigenvector. From this we conclude that all the components of $x$ are nonzero.

Now since $x$ is a real vector, we can define
$$V_1=\{ i\in [n] \ | \ x_i<0 \} $$
Then $V_1\ne \phi$ by (2.5) and $V_1\ne [n]$ by (2.5) and the fact that $k$ is even. This means that $V_1$ is a proper subset of the vertex set of $G$.

Finally, from (2.5) we see that for each edge $\{i_1,\cdots,i_k\} \in E(G)$, we have $x_{i_1}\cdots x_{i_k}<0$, which means that $|\{i_1,\cdots,i_k\}\cap V_1|$ is an odd number. This shows that $G$ is odd-bipartite, and thus completes the proof.
\qed

\vskip 1.88cm

\section{The Laplacian spectra and signless Laplacian spectra of the Cartesian products of $k$-uniform hypergraphs }

\vskip 0.18cm

In this section, we study the Laplacian spectra and signless Laplacian spectra of the Cartesian products of $k$-hypergraphs. We first show that the Cartesian product $G\Box H$ of two odd-bipartite $k$-hypergraphs $G$ and $H$ is still odd-bipartite. Then we use the direct product of tensors defined in [10] to obtain the expressions of the Laplacian tensors and the signless Laplacian tensors of the Cartesian products of hypergraphs. Using a useful relation $(\mathbb {A}\otimes \mathbb {B})(\mathbb {C}\otimes \mathbb {D})=(\mathbb {AC})\otimes (\mathbb {BD})$ between the direct product and general product of tensors, we are able to obtain some Laplacian eigenvalues of the Cartesian product $G\Box H$ from that of $G$ and $H$, and determine that the Laplacian spectral radius (and the largest Laplacian H-eigenvalue) of $G\Box H$ is the sum of the Laplacian spectral radii (and the largest Laplacian H-eigenvalues) of $G$ and $H$, when $G$ and $H$ are both connected and odd-bipartite.

\vskip 0.18cm

\noindent {\bf Definition 3.1} ([2], The Cartesian product of hypergraphs): Let $G$ and $H$ be two hypergraphs.
Define the Cartesian product $G\Box H$ of $G$ and $H$ as: $V(G\Box H)=V(G)\times V(H)$, and $\{(i_1,j_1),\cdots, (i_r,j_r)\}\in E(G\Box H)$ if and only if one of the following two conditions holds:

\vskip 0.1cm

\noindent (1). $i_1=\cdots =i_r$ and $\{j_1,\cdots ,j_r\}\in E(H)$.

\vskip 0.1cm

\noindent (2). $j_1=\cdots =j_r$ and $\{i_1,\cdots ,i_r\}\in E(G)$.

\vskip 0.28cm

It is easy to see that $G\Box H$ is k-uniform if both $G$ and $H$ are k-uniform, and is connected if both $G$ and $H$ are connected.

\vskip 0.18cm

\noindent {\bf Proposition 3.1}: Let $G$ and $H$ be two $k$-uniform odd bipartite hypergraphs with $k$ even. Then their Cartesian product $G\Box H$ is also odd bipartite.

\vskip 0.18cm

\noindent {\bf Proof.} Let $V_1(G)$ be the proper subset of $V(G)$ such that every edge of $G$ intersects $V_1(G)$ with exactly an odd number of vertices  (similar definition for $V_1(H)$). Let
$$V_1(G\Box H)=(V_1(G)\times V_1(H))\cup (\overline{V_1(G)}\times \overline{V_1(H)}).$$
We now show that every edge $e$ of $G\Box H$ intersects $V_1(G\Box H)$ with exactly an odd number of vertices.

\vskip 0.1cm

\noindent {\bf Case 1:} $e=\{(i,j_1),\cdots,(i,j_k)\}$ with $i\in V(G)$ and $\{j_1,\cdots,j_k\}\in E(H)$.

Without loss of generality, we may assume that $\{j_1,\cdots,j_k\}\cap V_1(H)=\{j_1,\cdots,j_r\}$, where $r$ is odd.

\vskip 0.1cm

\noindent {\bf Subcase 1.1:} If $i\in V_1(G)$. Then we have $e\cap V_1(G\Box H)= \{(i,j_1),\cdots,(i,j_r)\}$ whose cardinality $r$ is odd.

\vskip 0.1cm

\noindent {\bf Subcase 1.2:} If $i\notin V_1(G)$. Then we have $e\cap V_1(G\Box H)= \{(i,j_{r+1}),\cdots,(i,j_k)\}$ whose cardinality $k-r$ is also odd (since $k$ is even).

\vskip 0.1cm

\noindent {\bf Case 2:} $e=\{(i_1,j),\cdots,(i_k,j)\}$ with $j\in V(H)$ and $\{i_1,\cdots,i_k\}\in E(G)$. The proof of this case is similar to that of Case 1.

\vskip 0.1cm
\qed

\vskip 0.18cm

Now we consider the adjacency tensor, Laplacian tensor and signless Laplacian tensor of the Cartesian product $G\Box H$. First we recall the following concept of direct products of tensors defined in [10] which is a generalization of the direct products of matrices.

\vskip 0.18cm

\noindent {\bf Definition 3.2} [10]: Let $\mathbb {A}$ and $\mathbb {B}$ be two order $k$ tensors with dimension $n$ and $m$, respectively. Define the direct product $\mathbb {A}\otimes \mathbb {B}$ to be the following tensor of order $k$ and dimension $nm$ (the set of subscripts is taken as $[n]\times [m]$ in the lexicographic order):
$$(\mathbb {A}\otimes \mathbb {B})_{(i_1,j_1)(i_2,j_2)\cdots (i_k,j_k)}=a_{i_1i_2\cdots i_k}b_{j_1j_2\cdots j_k} $$

\vskip 0.2cm

The following relation between the direct product of tensors and the general product of tensors (defined in Section 1) can be found in [10] (This relation is also a generalization of a well-known similar relation for matrices).

\vskip 0.38cm

\noindent {\bf Proposition 3.2}: Let $\mathbb {A}$ and $\mathbb {B}$ be two order $k$ tensors with dimension $n$ and $m$, respectively. Let $\mathbb {C}$ and $\mathbb {D}$ be two order $r$ tensors with dimension $n$ and $m$, respectively. Then we have:
$$(\mathbb {A}\otimes \mathbb {B})(\mathbb {C}\otimes \mathbb {D})=(\mathbb {AC})\otimes (\mathbb {BD}).$$

\vskip 0.38cm

Using the direct product of tensors, the adjacency tensor, the Laplacian tensor and the signless Laplacian tensor of the Cartesian product of $k$-uniform hypergraphs can be obtained as in the following theorem.

\vskip 0.28cm

\noindent {\bf Theorem 3.1}: Let $\mathbb {A}(G)$ and $\mathbb {A}(H)$ be the adjacency tensors of a $k$-uniform hypergraph $G$ with $n$ vertices and a $k$-uniform hypergraph $H$ with $m$ vertices, respectively. Let $\mathbb {D}(G)$ and $\mathbb {D}(H)$ be the degree diagonal tensors of $G$ and $H$,  $\mathbb {L}(G)=\mathbb {D}(G)-\mathbb {A}(G)$ and $\mathbb {L}(H)=\mathbb {D}(H)-\mathbb {A}(H)$ be the Laplacian tensors of $G$ and $H$, $\mathbb {Q}(G)=\mathbb {D}(G)+\mathbb {A}(G)$ and $\mathbb {Q}(H)=\mathbb {D}(H)+\mathbb {A}(H)$ be the signless Laplacian tensors of $G$ and $H$. Let the ordering of the vertices of $G\Box H$ be taken as the lexicographic ordering of the elements of the set $V(G)\times V(H)$. Then we have:

\vskip 0.18cm

\noindent (1) ([10]). The adjacency tensor of $G\Box H$ is $\mathbb {A}(G\Box H)= \mathbb {A}(G)\otimes \mathbb {I}_m +  \mathbb {I}_n\otimes \mathbb {A}(H)$.

\vskip 0.18cm

\noindent (2). The degree diagonal tensor of $G\Box H$ is $\mathbb {D}(G\Box H)= \mathbb {D}(G)\otimes \mathbb {I}_m +  \mathbb {I}_n\otimes \mathbb {D}(H)$.

\vskip 0.18cm

\noindent (3). The Laplacian tensor of $G\Box H$ is $\mathbb {L}(G\Box H)= \mathbb {L}(G)\otimes \mathbb {I}_m +  \mathbb {I}_n\otimes \mathbb {L}(H)$.

\vskip 0.18cm

\noindent (4). The signless Laplacian tensor of $G\Box H$ is $\mathbb {Q}(G\Box H)= \mathbb {Q}(G)\otimes \mathbb {I}_m +  \mathbb {I}_n\otimes \mathbb {Q}(H)$.

\vskip 0.18cm

\noindent {\bf Proof. } (1).  Let $\mathbb {A}$, $\mathbb {B}$ and $\mathbb {C}$ be the adjacency tensors of the hypergraphs $G$, $H$ and $G\Box H$. Then by definition we can check that:
$$c_{(i_1,j_1)\cdots (i_k,j_k)}=\left \{\begin{array}{cc}
                               b_{j_1\cdots j_k}&   \mbox{if }  \ i_1=i_2=\cdots =i_k  \\
                                a_{i_1\cdots i_k} & \mbox{ if } \  j_1=j_2=\cdots =j_k  \\
                               0 &   \mbox{otherwise}\\
                            \end{array}
                          \right.  \eqno {(3.1)}$$
Notice that all the diagonal entries of $\mathbb {A}$ and $\mathbb {B}$ are zero, so it follows from (3.1) that
$$c_{(i_1,j_1)\cdots (i_k,j_k)}=a_{i_1\cdots i_k}\delta _{j_1\cdots j_k}+\delta_{i_1\cdots i_k}b_{j_1\cdots j_k}$$
It is also easy to see that
$$(\mathbb {A}\otimes \mathbb {I}_m +  \mathbb {I}_n\otimes \mathbb {B})_{(i_1,j_1)\cdots (i_k,j_k)}=a_{i_1\cdots i_k}\delta _{j_1\cdots j_k}+\delta_{i_1\cdots i_k}b_{j_1\cdots j_k}  \eqno {(3.2)}$$
So we have $\mathbb {C}=\mathbb {A}\otimes \mathbb {I}_m +  \mathbb {I}_n\otimes \mathbb {B}$.

\vskip 0.18cm

\noindent (2). For two diagonal tensors, we only need to compare their diagonal entries. We have

$$\mathbb {D}(G\Box H)_{(i,j)\cdots (i,j)}=d_{G\Box H} ((i,j))=d_G(i)+d_H(j)   \eqno {(3.3)}$$
On the other hand, we have
$$ (\mathbb {D}(G)\otimes \mathbb {I}_m +  \mathbb {I}_n\otimes \mathbb {D}(H) )_{(i,j)\cdots (i,j)} = \mathbb {D}(G)_{i\cdots i}+\mathbb {D}(H)_{j\cdots j} = d_G(i)+d_H(j)  \eqno {(3.4)} $$
By (3.3) and (3.4) we have
$$\mathbb {D}(G\Box H)_{(i,j)\cdots (i,j)}=  (\mathbb {D}(G)\otimes \mathbb {I}_m +  \mathbb {I}_n\otimes \mathbb {D}(H) )_{(i,j)\cdots (i,j)}  $$
Thus we obtain $\mathbb {D}(G\Box H)= \mathbb {D}(G)\otimes \mathbb {I}_m +  \mathbb {I}_n\otimes \mathbb {D}(H)$ (since both sides are diagonal tensors).

\vskip 0.18cm

\noindent (3). From (1) and (2), we can easily obtain:

$$\begin{aligned}
\displaystyle &  \mathbb {L}(G\Box H)= \mathbb {D}(G\Box H) - \mathbb {A}(G\Box H)   \\
&= (\mathbb {D}(G)\otimes \mathbb {I}_m +  \mathbb {I}_n\otimes \mathbb {D}(H))-(\mathbb {A}(G)\otimes \mathbb {I}_m +  \mathbb {I}_n\otimes \mathbb {A}(H))   \\
&= (\mathbb {D}(G)-\mathbb {A}(G))\otimes \mathbb {I}_m +  \mathbb {I}_n\otimes (\mathbb {D}(H)-\mathbb {A}(H))        \\
&= \mathbb {L}(G)\otimes \mathbb {I}_m +  \mathbb {I}_n\otimes \mathbb {L}(H)\\
\end{aligned}
$$

\vskip 0.1cm

\noindent The proof of (4) is the same as that of (3).
\qed

\vskip 0.18cm

By using Proposition 3.2, we can obtain the following relation between the eigenvalue-eigenvectors of the tensors $\mathbb {A}$ and $\mathbb {B}$ and that of $\mathbb {A}\otimes \mathbb {I}_m +  \mathbb {I}_n\otimes \mathbb {B}$.

\vskip 0.38cm

\noindent {\bf Theorem 3.2}: Let $\mathbb {A}$ and $\mathbb {B}$ be two order $k$ tensors with dimension $n$ and $m$, respectively. Suppose that we have $\mathbb {A}u=\lambda u^{[k-1]}$, and $\mathbb {B}v=\mu v^{[k-1]}$, and write $w=u\otimes v$. Then we have:  $(\mathbb {A}\otimes \mathbb {I}_m +  \mathbb {I}_n\otimes \mathbb {B})w=(\lambda +\mu)w^{[k-1]}.$

\vskip 0.1cm

\noindent {\bf Proof. }  We have by Proposition 3.2 and Theorem 3.1 that:

\vskip 0.1cm

\noindent
$$\begin{aligned}
\displaystyle &   (\mathbb {A}\otimes \mathbb {I}_m +  \mathbb {I}_n\otimes \mathbb {B})w=(\mathbb {A}\otimes \mathbb {I}_m +  \mathbb {I}_n\otimes \mathbb {B})(u\otimes v)= (\mathbb {A}u)\otimes (\mathbb {I}_m v)+ ( \mathbb {I}_n u)\otimes (\mathbb {B}v)\\
&=(\lambda u^{[k-1]})\otimes (v^{[k-1]}) + ( u^{[k-1]})\otimes (\mu v^{[k-1]})  \\
&=(\lambda +\mu)(u^{[k-1]}\otimes v^{[k-1]})=(\lambda +\mu)(u\otimes v)^{[k-1]}= (\lambda +\mu)w^{[k-1]}
\qquad \qquad \qquad \qquad \qed \\
\end{aligned}
$$

From Theorem 3.2, we can obtain the following results about the Laplacian spectra (and the signless Laplacian spectra) of the hypergraphs $G\Box H$.

\vskip 0.38cm

\noindent {\bf Corollary 3.1}: Let $G$ and $H$ be two $k$-uniform hypergraphs. Let $\lambda $ be a Laplacian eigenvalue (or signless Laplacian eigenvalue) of $G$ with eigenvector $u$, and $\mu $ be a Laplacian eigenvalue (or signless Laplacian eigenvalue) of $H$ with eigenvector $v$, respectively. Then  $ \ \lambda +\mu$ is a Laplacian eigenvalue (or signless Laplacian eigenvalue) of $G\Box H$ with eigenvector $u\otimes v$.

\vskip 0.38cm

\noindent {\bf Theorem 3.3}: Let $\mathbb {A}\ge O$ and $\mathbb {B}\ge O$ be two order $k$ nonnegative tensors with dimension $n$ and $m$, respectively. Then we have:

\vskip 0.18cm

\noindent (1). $\rho (\mathbb {A}\otimes \mathbb {I}_m +  \mathbb {I}_n\otimes \mathbb {B})=\rho (\mathbb {A})+\rho (\mathbb {B})$.

\vskip 0.18cm

\noindent (2). $\lambda (\mathbb {A}\otimes \mathbb {I}_m +  \mathbb {I}_n\otimes \mathbb {B})=\lambda (\mathbb {A})+\lambda (\mathbb {B})$.

\vskip 0.18cm

\noindent {\bf Proof. } (1). We consider the following two cases.

\vskip 0.18cm

\noindent {\bf Case 1}. Both $\mathbb {A}$ and $\mathbb {B}$ are positive tensors.

By Lemma 1.1, let $u$ (and $v$) be the positive eigenvector of the tensor $\mathbb {A}$ (and $\mathbb {B}$) corresponding to the eigenvalue $\rho (\mathbb {A})$ (and $\rho (\mathbb {B})$), respectively (since $\mathbb {A}$ and $\mathbb {B}$ are both positive tensors). Then by Theorem 3.2, $\rho (\mathbb {A})+ \rho (\mathbb {B})$ is an eigenvalue of $\mathbb {A}\otimes \mathbb {I}_m +  \mathbb {I}_n\otimes \mathbb {B}$ with a positive eigenvector $u\otimes v$. Thus from (2) of Lemma 1.1 we see that $\rho (\mathbb {A})+\rho (\mathbb {B})$ must be the spectral radius of $\mathbb {A}\otimes \mathbb {I}_m +  \mathbb {I}_n\otimes \mathbb {B}$. This proves Case 1.

\vskip 0.18cm

\noindent {\bf Case 2}. The general case (when $\mathbb {A}$ and $\mathbb {B}$ are both nonnegative).

Take $\epsilon >0$ and $\mathbb {A}_{\epsilon }=\mathbb {A}+\epsilon \mathbb {J}_1$ and $\mathbb {B}_{\epsilon }=\mathbb {B}+\epsilon \mathbb {J}_2$, where $\mathbb {J}_1$ and $\mathbb {J}_2$ are order $k$ tensors with all entries 1 with dimension $n$ and $m$, respectively. Then $\mathbb {A}_{\epsilon }$ and $\mathbb {B}_{\epsilon }$ are both positive tensors. So by Case 1 we have
$$\rho (\mathbb {A}_{\epsilon }\otimes \mathbb {I}_m +  \mathbb {I}_n\otimes \mathbb {B}_{\epsilon })=\rho (\mathbb {A}_{\epsilon })+\rho (\mathbb {B}_{\epsilon })$$
Take the limit $\epsilon \rightarrow 0$ on both sides of the above equation (since the maximal absolute value of the roots of a complex polynomial is a continuous function on the coefficients of the polynomial), we obtain the desired result.

\vskip 0.18cm

\noindent (2). By Lemma 1.1 we know that for any nonnegative tensor $\mathbb {T}$, we have $\lambda (\mathbb {T})=\rho (\mathbb {T})$. Thus (2) follows directly from (1).
\qed

\vskip 0.18cm

The following theorem shows that the Laplacian spectral radius (and the largest Laplacian H-eigenvalue) of $G\Box H$ is the sum of the Laplacian spectral radii (and the largest Laplacian H-eigenvalues) of $G$ and $H$, when $G$ and $H$ are both connected odd-bipartite.

\vskip 0.38cm

\noindent {\bf Theorem 3.4}: Let $G$ and $H$ be two $k$-uniform hypergraphs with $n$ vertices and $m$ vertices, respectively. Let $\mathbb {A}(G)$ and $\mathbb {A}(H)$ be the adjacency tensors of  $G$ and $H$, $\mathbb {L}(G)$ and $\mathbb {L}(H)$ be the Laplacian tensors of $G$ and $H$, $\mathbb {Q}(G)$ and $\mathbb {Q}(H)$ be the signless Laplacian tensors of $G$ and $H$, respectively. Then we have:

\vskip 0.18cm

\noindent (1). $\rho (\mathbb {A}(G\Box H))=\rho (\mathbb {A}(G))+\rho (\mathbb {A}(H))$, and $\lambda (\mathbb {A}(G\Box H))=\lambda (\mathbb {A}(G))+\lambda (\mathbb {A}(H))$.

\vskip 0.18cm

\noindent (2). $\rho (\mathbb {Q}(G\Box H))=\rho (\mathbb {Q}(G))+\rho (\mathbb {Q}(H))$, and $\lambda (\mathbb {Q}(G\Box H))=\lambda (\mathbb {Q}(G))+\lambda (\mathbb {Q}(H))$.

\vskip 0.18cm

\noindent (3). If we further assume that $k$ is even, and $G$ and $H$ are both connected and odd-bipartite, then we have: $\rho (\mathbb {L}(G\Box H))=\rho (\mathbb {L}(G))+\rho (\mathbb {L}(H))$, and $\lambda (\mathbb {L}(G\Box H))=\lambda (\mathbb {L}(G))+\lambda (\mathbb {L}(H))$.

\vskip 0.2cm

\noindent {\bf Proof. } (1)+(2). Since $\mathbb {A}(G)$, $\mathbb {A}(H)$, $\mathbb {Q}(G)$ and $\mathbb {Q}(H)$ are all nonnegative, (1) and (2) follow directly from Theorem 3.3 and Theorem 3.1.

\vskip 0.18cm

\noindent (3). Since $G$ and $H$ are both connected and odd-bipartite, we see from Proposition 3.1 that $G\Box H$ is also connected and odd-bipartite. Thus by Theorem 2.2 we have:
$$\rho (\mathbb {L}(G\Box H))= \rho (\mathbb {Q}(G\Box H)),  \   \rho (\mathbb {L}(G))=\rho (\mathbb {Q}(G)), \ \ \rho (\mathbb {L}(H))=\rho (\mathbb {Q}(H)) \eqno {(3.5)}$$
and
$$\lambda (\mathbb {L}(G\Box H))= \lambda (\mathbb {Q}(G\Box H)),  \   \lambda (\mathbb {L}(G))=\lambda (\mathbb {Q}(G)), \ \lambda (\mathbb {L}(H))=\lambda (\mathbb {Q}(H)) \eqno {(3.6)}$$
From (3.5) and the result (2), we obtain
$$\rho (\mathbb {L}(G\Box H))=\rho (\mathbb {Q}(G\Box H))=\rho (\mathbb {Q}(G))+\rho (\mathbb {Q}(H))=\rho (\mathbb {L}(G))+\rho (\mathbb {L}(H)).$$
From (3.6) and the result (2), we obtain
$$\lambda (\mathbb {L}(G\Box H))=\lambda (\mathbb {Q}(G\Box H))=\lambda (\mathbb {Q}(G))+\lambda (\mathbb {Q}(H))=\lambda (\mathbb {L}(G))+\lambda (\mathbb {L}(H)).$$
\qed

\vskip 0.88cm

\section{Final Remarks}

\vskip 0.18cm

From Theorem 2.2, we can see that if a $k$-hypergraph $G$ is connected, then we have the following implication:
$$Hspec(\mathbb {L})= Hspec (\mathbb {Q}) \Longrightarrow Spec(\mathbb {L})= Spec (\mathbb {Q})$$
But we do not know whether the reverse implication
is true or not (since H-eigenvalues need to have real eigenvectors).

\vskip 0.18cm

On the other hand, if the reverse implication is true, then the condition (2) of Theorem 2.2 could be simply replaced by $Spec(\mathbb {L})= Spec (\mathbb {Q})$, and this would enable us to provide more connection between Theorem 2.2 and Theorem 2.4.

\vskip 1.8cm
\noindent{\bf References}

\bibliographystyle{acm}
%\bibliography{reference}

\end{document}